\documentclass[11pt,a4paper]{article}
\usepackage[pdftex]{graphicx}
\usepackage{amsmath,amssymb,amsfonts,amsthm}
\numberwithin{equation}{section}
\usepackage{indentfirst}
\usepackage{enumitem} 
\usepackage[margin=1.3in]{geometry}
\usepackage{fancyhdr}

\usepackage[colorlinks=true,linkcolor=magenta,citecolor=blue,urlcolor=cyan]{hyperref}
\usepackage{titlesec}
\usepackage{etoolbox}
\usepackage{graphicx}  
\usepackage{changepage}
\usepackage{float}
\theoremstyle{plain}
\newtheorem{theorem}{Theorem}[section]
\newtheorem{lemma}[theorem]{Lemma}

\newtheorem{definition}[theorem]{Definition}

\newcommand{\leg}[2]{\left(\frac{#1}{#2}\right)}
\usepackage{comment}

\newcommand{\Z}{\mathbb{Z}}
\newcommand{\RD}[2]{\mathrm{RD}_{l,t}^{#1}(#2)}
\makeatletter
\renewcommand{\maketitle}{
	\begin{center}
		{\Large\bfseries{\@title}\par}
		\vskip 1em
		{\normalsize
			\lineskip .5em
			\begin{tabular}[t]{c}
				\@author
			\end{tabular}\par}
		\vskip 1.5em
	\end{center}
}
\makeatother

\usepackage{color}
\usepackage{xcolor}
\usepackage[normalem]{ulem} 
\usepackage{soul}

\usepackage{graphicx} 
\usepackage{xstring}  

\usepackage{xstring}
\usepackage{graphicx} 

\usepackage{marvosym}

\begin{document}
	\begin{center}{\bf \large
			{New Congruences for Simultaneously $l$-Regular and $t$-Distinct Partition in $k$-Colors \\[0.15cm] }  }\vspace{0.2cm}
	\end{center}
	\begin{center}
		\footnotemark[1]
		\bf Anjelin Mariya Johnson, \footnotemark[2]
		\bf S. N. Fathima  

	\end{center}

	\vspace{1 cm}
	\begin{center}
		\begin{minipage}{0.85\textwidth}  
			
			\noindent{\bf Abstract:} 
			   Let $RD_{\ell,t}^k(n)$ enumerates the number of $k$-colored simultaneously $l$-regular and  $t$-distinct partitions. In this article, we establish some new infinite families of congruences modulo prime $p$, $p\geq5$ and modulo powers of 5 for $RD_{\ell,t}^k(n)$. For instance, 
			   for each $n\geq0$, $0\leq\alpha\leq{2k-2}$ and $k\geq1$,
			   \begin{align*}
			   	\sum_{n=0}^{\infty}RD_{(5^{2k-1},5^{2k-1})}^2\big(5^{2k+\alpha-1}n+\frac{9\cdot5^{2k+\alpha-1}-2\cdot5^{2k-1}+1}{12}\big)q^n\equiv& 0\pmod{5^{k+\left[\frac{\alpha}{2}\right]}}
			   \end{align*}
			\vspace{.125cm}
			
			\vspace{.05cm}
			
			\end{minipage}
		\end{center}

\section {Introduction}
For positive integer $n$ and a complex number $q$ with $|q|<1$, $E_r:=\prod_{n\geq1}(1-q^{nr})$. 
A partition of a nonnegative integer $n$ is defined by a sequence of positive integers $(\lambda_1, \lambda_2,\ldots,\lambda_r)$ with $\lambda_1\geq\lambda_2\geq\ldots\geq\lambda_r$ such that $n=\lambda_1+\lambda_2+\ldots+\lambda_r$. The integers $\lambda_r$ are called parts of the partition. Closely related to $q$-series products is Rogers Ramanujan continued fraction $R(q),$\cite{power}:
\begin{align*}
	R(q):=\frac{(q;q^5)_\infty(q^4;q^5)_\infty}{(q^2;q^5)_\infty(q^3;q^5)_\infty}.
\end{align*}
which satisfies the following 5 dissection for $E_1$ involving $R(q)$
\begin{align}
	E_1=E_{25}\left(\frac{1}{R(q^5)}-q-q^2R(q^5)\right).\label{r(q)}
\end{align}
	Keith \cite{keith}, introduced  simultaneously $l$-regular and $t$-distinct partition function, $RD_{\ell,t}(n)$.
A partition of positive integer $n$ is said to be simultaneously $l$-regular and $t$-distinct if none of the parts are divisible by $l$ and parts appear fewer than $t$ times. The generating function for $RD_{\ell,t}$ is given by
\begin{align*}
	\sum_{n=0}^{\infty} RD_{\ell,t}(n)q^n=
	\frac{E_{\ell}E_t}{E_1E_{\ell t}}.
\end{align*}
For example, four partitions of $RD_{3,4}(5)$ are: $(5),\;(4,1),\;(2,2,1),\;(2,1,1,1).$\\
Over the last several years, the partition function $RD_{\ell,t}(n)$ has received considerable attention owing to its rich arithmetic properties. For more information on this partition function $RD_{\ell,t}(n)$, we refer the readers to \cite{keith,keith2,Ahmia 2}. 
Recently, the authors \cite{Anjelin} considered the family of functions $RD_{\ell,t}^k$ which counts the $k$-colored simultaneously $l$-regular and  $t$-distinct partitions. The generating function for $RD_{\ell,t}^k(n)$ is given by
\begin{align}\label{Gf}
	\sum_{n=0}^{\infty} RD_{\ell,t}^k(n)q^n=\left(\frac{E_{\ell}E_t}{E_1E_{\ell t}}\right)^k.
\end{align}
For example, $RD_{3,4}^3(2)=9$, with the relevant partitions being
\begin{align*}
	(2_1),(2_2),(2_3),(1_1,1_1),(1_2,1_2),(1_3,1_3),(1_1,1_2),(1_1,1_3),(1_2,1_3).
\end{align*}

Our goal in this brief work is to specifically revisit the function $RD_{\ell,t}^k(n)$ to establish new infinite families of congruences modulo prime $p$, $p\geq5$ and modulo powers of 5.
We now state our results.\\
The following theorem gives infinite families of congruences modulo $p$, $p\geq5$.
\begin{theorem}\label{thm:main}
	Let $p\ge5$ be prime, let $l,t$ be positive integers with $p\mid l$ and
	$p\mid t$, and let $k\ge1$. Then for all $n\ge0$ and $0\le r\le p-1$,
	\begin{align*}
		RD_{l,t}^{k}({pn+r})\equiv 0 \pmod p
	\end{align*}
	\begin{align*}
		\quad\text{whenever}\quad
		\begin{cases}
			k\equiv 0 \pmod p, & r\not\equiv 0 \pmod p,\\[6pt]
			k\equiv -1 \pmod p, & \leg{24r+1}{p}=-1,\\[6pt]
			k\equiv -3 \pmod p, & \leg{8r+1}{p}\in\{-1,0\},\\[6pt]
			k\equiv -4 \pmod p,\ p\equiv 2 \pmod 3, & 6r+1\equiv 0 \pmod p.
		\end{cases}
	\end{align*}
	Moreover, if $k=pm$, then
	\begin{align*}
		RD_{l,t}^{pm}({pn})\equiv RD_{l,t}^{m}({n})\pmod{p}.
	\end{align*}
\end{theorem}
 \noindent The following theorem gives infinite families of congruences modulo powers of 5.
\begin{theorem}\label{T1}
For each $n\geq0$, $0\leq\alpha\leq{2k-2}$ and $k\geq1$,
\begin{align}
	\sum_{n=0}^{\infty}RD_{(5^{2k-1},5^{2k-1})}^2\big(5^{2k+\alpha-1}n+\frac{9\cdot5^{2k+\alpha-1}-2\cdot5^{2k-1}+1}{12}\big)q^n\equiv& 0\pmod{5^{k+\left[\frac{\alpha}{2}\right]}}\label{t11},\\
	\sum_{n=0}^{\infty}RD_{(5^{2k-1},5^{2k-1})}^2\big(5^{2k+\alpha}n+\frac{33\cdot5^{2k+\alpha-1}-2\cdot5^{2k-1}+1}{12}\big)q^n\equiv& 0\pmod{5^{k+\left[\frac{\alpha}{2}\right]+1}}\label{t12},\\
	\sum_{n=0}^{\infty}RD_{(5^{2k-1},5^{2k-1})}^2\big(5^{2k+\alpha}n+\frac{57\cdot5^{2k+\alpha-1}-2\cdot5^{2k-1}+1}{12}\big)q^n\equiv& 0\pmod{5^{k+\left[\frac{\alpha}{2}\right]+1}}\label{t14}.
		\end{align}
	
\end{theorem}
\begin{theorem}\label{t2}
	For each $n\geq0$, $\alpha\geq{0}$ and $k\geq1$,
	\begin{align}
		\sum_{n=0}^{\infty}RD_{(5^{2k-1},5^{2k-1})}^2\big(5^{4k+2\alpha-2}n+\frac{8\cdot5^{4k+2\alpha-2}+({5^{2k-1}-1})^2}{12}\big)q^n\equiv0& \pmod{5^{2k+\alpha}}\label{t21},\\
		\sum_{n=0}^{\infty}RD_{(5^{2k-1},5^{2k-1})}^2\big(5^{4k+2\alpha-1}n+\frac{4\cdot5^{4k+2\alpha-1}+({5^{2k-1}-1})^2}{12}\big)q^n\equiv0& \pmod{5^{2k+\alpha}}\label{t22},\\
		\sum_{n=0}^{\infty}RD_{(5^{2k-1},5^{2k-1})}^2\big(5^{4k+2\alpha-1}n+\frac{44\cdot5^{4k+2\alpha-2}+({5^{2k-1}-1})^2}{12}\big)q^n\equiv0& \pmod{5^{2k+\alpha+1}}\label{t23},\\
		\sum_{n=0}^{\infty}RD_{(5^{2k-1},5^{2k-1})}^2\big(5^{4k+2\alpha-1}n+\frac{56\cdot5^{4k+2\alpha-2}+({5^{2k-1}-1})^2}{12}\big)q^n\equiv0& \pmod{5^{2k+\alpha+1}}\label{t24},\\
		\sum_{n=0}^{\infty}RD_{(5^{2k-1},5^{2k-1})}^2\big(5^{4k+2\alpha}n+\frac{28\cdot5^{4k+2\alpha-1}+({5^{2k-1}-1})^2}{12}\big)q^n\equiv0& \pmod{5^{2k+\alpha+1}}\label{t25},\\
		\sum_{n=0}^{\infty}RD_{(5^{2k-1},5^{2k-1})}^2\big(5^{4k+2\alpha}n+\frac{52\cdot5^{4k+2\alpha-1}+({5^{2k-1}-1})^2}{12}\big)q^n\equiv0& \pmod{5^{2k+\alpha+1}} \label{t26}.
	\end{align}
\end{theorem}

	\noindent Our article is structured as follows: In Section 2 we collect all of the tools necessary to complete the proofs of our results. Section 3 will focus on generating functions for simultaneously $l$-regular and  $t$-distinct partitions in 2 colors. Finally in section 4 we provide the proof of Theorem 1.1 and slightly different proof of Theorems 1.2 and 1.3 using the techniques of $H$ operator first employed by Hirschhorn.
\section{Preliminaries} 
In this section, we list certain definitions and lemmas that plays a vital role in proving our results.
\begin{definition}
	Let $p$ be an odd prime and $a \in \mathbb{Z}$. The \emph{Legendre symbol}
	$\left(\frac{a}{p}\right)$ is defined by
	\[
	\left(\frac{a}{p}\right) :=
	\begin{cases}
		\phantom{-}1 & \text{if } a \text{ is a quadratic residue mod } p \text{ and } p \nmid a,\\
		\phantom{-}0 & \text{if } p \mid a,\\
		-1 & \text{if } a \text{ is a quadratic non-residue mod } p.
	\end{cases}
	\]
\end{definition}
Using Binomial theorem, we obtain a congruence result which is used frequently in the sequel without explicitly mentioning.
\begin{align}
	E_m^{p^k} \equiv E_{mp}^{p^{k-1}} \pmod{p^k}.\label{binomial}
\end{align}
We next recall a few classical results from the theory of $q$-series.
\begin{lemma}[Euler's Pentagonal number Theorem {\cite[(1.6.1)]{power}}]
\begin{align}
	E_1=\sum_{k=-\infty}^{\infty}(-1)^kq^\frac{k(3k+1)}{2}.\label{pentagonal}
\end{align}
\end{lemma}
\begin{lemma}[Jacobi's identity {\cite[(1.7.1)]{power}}]
\begin{align}
	E_1^3 = \sum_{n \geq 0} (-1)^n (2n+1) q^{n(n+1)/2}. \label{jacobi}
\end{align}
\end{lemma}
Thanks to (\ref{r(q)}) we define
\begin{align*}
	\zeta = \frac{E_1}{qE_{25}}, \qquad T = \frac{E_5^6}{q^5 E_{25}^6}.
\end{align*}

Following Hirschhorn \cite{power, Hunt}, we introduce a ``huffing'' operator $H$, which operates on a series by picking out those terms of the form $q^{5n}$, and huffing the rest away. That is,
\[
H\left( \sum_{n=0}^{\infty} a_n q^n \right) = \sum_{n=0}^{\infty} a_{5n} q^{5n}.
\]
As in Hirschhorn \cite{power}, we define an infinite matrix $\{m(i,j)\}_{i,j\geq 1}$ by
\[
\begin{pmatrix}
	5 & 0 & 0 & 0 & 0 & 0 & \cdots\\
	2\times 5 & 5^3 & 0 & 0 & 0 & 0 & \cdots\\
	9 & 3\times 5^3 & 5^5 & 0 & 0 & 0 & \cdots\\
	4 & 22\times 5^2 & 4\times 5^5 & 5^7 & 0 & 0 & \cdots\\
	1 & 4\times 5^3 & 8\times 5^5 & 5^8 & 5^9 & 0 & \cdots\\
	\vdots & \vdots & \vdots & \vdots & \vdots & \vdots & \ddots
\end{pmatrix}
\]
and for $i\geq 6$, $m(i,1)=0$, and for $j\geq 2$,
\begin{align}
	m(i,j) &= 25m(i-1,j-1) + 25m(i-2,j-1) + 15m(i-3,j-1) \notag\\
	&\qquad + 5m(i-4,j-1) + m(i-5,j-1). \label{eq:2.4}
\end{align}
By induction, it follows immediately that

\begin{lemma}[\cite{power}] \label{lem:2.1}
	We have
	\begin{itemize}
		\item $m(i,j)=0$ for $j>i$.
		\item $m(i,j)=0$ for $i>5j$.
	\end{itemize}
\end{lemma}

The following lemma is important for our proof.

\begin{lemma}[\cite{power}, Eq.~(6.4.9)] \label{lem:2.2}
	For $i\geq 1$, we have
	\[
	H\left(\frac{1}{\zeta^i}\right) = \sum_{j=1}^{\infty} \frac{m(i,j)}{T^j} = \sum_{j=1}^{i} \frac{m(i,j)}{T^j}.
	\]
\end{lemma}

\begin{lemma}
	For integers $i \geq 1$, we have
\begin{align}
		H\!\left(q^{i-4}\frac{E_{5}^{6i-2}}{E_{1}^{6i-4}}\right)=&
	\sum_{j=1}^{5i-3} a_{i,j}\,\label{L1}
	q^{5j-5}\frac{E_{25}^{6j-2}}{E_{5}^{6j-4}},\\
		H\!\left(q^{i-4}{\big(\frac{E_{5}}{E_{1}}\big)}^{6i-4}\right)=&
	\sum_{j=1}^{5i-3} a_{i,j}\, \label{L2}
	q^{5j-5}{(\frac{E_{25}}{E_{5}})}^{6j-2},\\
	H\!\left(q^{i-2}{\big(\frac{E_{5}}{E_{1}}\big)}^{6i-2}\right)=&
	\sum_{j=1}^{5i-1} b_{i,j}\, \label{L3}
	q^{5j-5}{(\frac{E_{25}}{E_{5}})}^{6j-4}. 
\end{align}
\end{lemma}
\begin{proof}
	Note for $1\leq j \leq i-1$, $m_{6i-4,j}=0$, therefore, we have
	\noindent
	\begin{align*}\label{a15}
		H\!\left(\frac{1}{\zeta^{6i-4}}\right)
		= \sum_{j= i}^{6i-4} \frac{m_{6i-4,\,j}}{T^{j}}
		= \sum_{j= 1}^{5i-3} \frac{m_{6i-4,\,i+j-1}}{T^{i+j-1}},
	\end{align*}
	which can be rewritten as
	\begin{equation*}
		H\!\left(q^{6i-4}\frac{E_{25}^{6i-4}}{E_{1}^{6i-4}}\right)
		=
		\sum_{j=1}^{5i-3}
		m_{6i-4,i+j-1}\,
		q^{5i+5j-5}
		\frac{E_{25}^{6i+6j-6}}{E_{5}^{6i+6j-6}},
	\end{equation*}
	that gives
	\begin{equation*}
	H\!\left(q^{i-4}\frac{E_{5}^{6i-2}}{E_{1}^{6i-4}}\right)=
	\sum_{j=1}^{5i-3} a_{i,j}\,
	q^{5j-5}\frac{E_{25}^{6j-2}}{E_{5}^{6j-4}}
	\end{equation*}
	where
	\[
	a_{i,j}=m_{6i-4,i+j-1},
	\]
	this completes the proof of (\ref{L1}). The identity (\ref{L2}) follows from identity (\ref{L1}).
	Similarly, for $1\leq j \leq i-1$, $m_{6i-2,j}=0$, then
	\begin{align*}
		H\!\left(\frac{1}{\zeta^{6i-2}}\right)
		= \sum_{j= i}^{6i-2} \frac{m_{6i-2,\,j}}{T^{j}}
		= \sum_{j= 1}^{5i-1} \frac{m_{6i-2,\,i+j-1}}{T^{i+j-1}},
	\end{align*}
	i.e,
	\begin{equation*}
		H\!\left(q^{6i-2}\frac{E_{25}^{6i-2}}{E_{1}^{6i-2}}\right)
		=
		\sum_{j=1}^{5i-1}
		m_{6i-2,i+j-1}\,
		q^{5i+5j-5}
		\frac{E_{25}^{6i+6j-6}}{E_{5}^{6i+6j-6}},
	\end{equation*}
	that gives
	\begin{equation*}
		H\!\left(q^{i-2}{(\frac{E_{5}}{E_{1}})}^{6i-2}\right)=
	\sum_{j=1}^{5i-1} b_{i,j}\,
	q^{5j-5}{(\frac{E_{25}}{E_{5}})}^{6j-4}. 
	\end{equation*}
	where
	\[
	b_{i,j}=m_{6i-2,i+j-1},
	\]
	this completes the proof of (\ref{L3}). 
\end{proof}
\section{Generating Functions}
In this section, we establish the generating functions for $RD_{\ell,t}^2$ within specific arithmetic progressions. We will heavily rely on this generating function in next section.
\begin{theorem}\label{gf1}
	For each $n\geq0$, $0\leq\alpha\leq{2k-2}$ and $k\geq1$,
	\begin{align}
		\sum_{n=0}^{\infty}RD_{(5^{2k-1},5^{2k-1})}^2\big(5^{2k+\alpha-1}n+\frac{9\cdot5^{2k+\alpha-1}-2\cdot5^{2k-1}+1}{12}\big)q^n \nonumber\\=\frac{1}{({E_{5^{2k-\alpha-1}}})^2}\sum_{i=1}^{\infty}y_{\alpha+1,i}^{2k-1}q^{i-1}\frac{E_5^{6i-2}}{E_1^{6i-4}}.
			\end{align}
			where the coefficient entries are defined as follows:
		\begin{align*}
				Y_1^{2k-1}=&X_{2k-1}
			\end{align*}
				and for any $1\leq\alpha\leq{2k-2}$,
			\begin{align}\label{Ydef}
			Y_{\alpha+1}^{2k-1}=&Y_\alpha^{2k-1}\cdot A,
		\qquad	A=(a_{i,j})=m_{6i-4,i+j-1}.
		\end{align}
		The coefficient vectors $X_1=(10,125,0,0,\ldots)$,
		and for any $s\geq1,$ we have\\
	\begin{align*}
			X_{s+1}=\begin{cases}
			 X_s \cdot P,& \text{if $s$ is odd}\\
			 X_s \cdot Q,& \text{if $s$ is even}
		\end{cases}
	\end{align*}
	where,
	\begin{align*}
			P:=&(p_{i,j})=m_{6i,i+j}\\
		Q:=&(q_{i,j})=m_{6i+2,i+j}
	\end{align*}
\end{theorem}
\begin{proof}
	It follows from (\ref{Gf}) that
	\begin{align}\label{3.1.1}
		\sum_{n=0}^{\infty} RD_{(5^{2k-1},5^{2k-1})}^2(n)\,q^n=\frac{{E_{5^{2k-1}}}^4}{{E_1}^2{E_{5^{4k-2}}}^2}
		=\frac{{E_{5^{2k-1}}}^4}{{E_{5^{4k-2}}}^2} \sum_{n=0}^{\infty} p_{-2}(n)\,q^n,
	\end{align}
	where $p_{-2}(n)$ is defined by
	\[
	\sum_{n=0}^{\infty} p_{-2}(n)\,q^n = \frac{1}{{E_1}^2}.
	\]
	Extracting the terms involving $q^{5^{2k-1}n+(7\times 5^{2k-1}+1)/12}$ in (\ref{3.1.1}) and upon simplification, yields
	\begin{align}
			\sum_{n=0}^{\infty} RD_{(5^{2k-1},5^{2k-1})}^2\!\left(5^{2k-1}n+\frac{7\times 5^{2k-1}+1}{12}\right)q^n\nonumber\\
			= \frac{{E_1}^4}{{E_{5^{2k-1}}}^2} \sum_{n=0}^{\infty} p_{-2}\!\left(5^{2k-1}n+\frac{7\times 5^{2k-1}+1}{12}\right)q^n. \label{after}
		\end{align}
	Thanks to Tang \cite[Theorem 4.1]{tang}, we have
	\begin{align}
		\sum_{n=0}^{\infty} p_{-2}\!\left(5^{2k-1}n+\frac{7\times 5^{2k-1}+1}{12}\right)q^n
		= \sum_{j=1}^{\infty} x_{2k-1,j}\, q^{j-1}\,\frac{E_5^{6j-2}}{E_1^{6j}},\label{tang}
	\end{align}
Employing (\ref{tang}) in (\ref{after}) yields
	\[
	\sum_{n=0}^{\infty} RD_{(5^{2k-1},5^{2k-1})}^2\!\left(5^{2k-1}n+\frac{7\times 5^{2k-1}+1}{12}\right)q^n
	= \frac{1}{{E_{5^{2k-1}}}^2}\sum_{j=1}^{\infty} x_{2k-1,j}\, q^{j-1}\,\frac{E_5^{6j-2}}{E_1^{6j-4}}.
	\]
	This implies that $\alpha=0$ case of Theorem \ref{gf1} holds.
	
	Now for some $\alpha$ with $0\leq\alpha\leq2k-3$ , we assume that Theorem \ref{gf1} holds, specifically
	\begin{align*}
		\sum_{n=0}^{\infty}RD_{(5^{2k-1},5^{2k-1})}^2\big(5^{2k+\alpha-1}n+\frac{9\cdot5^{2k+\alpha-1}-2\cdot5^{2k-1}+1}{12}\big)q^{n-3} \nonumber\\=\frac{1}{({E_{5^{2k-\alpha-1}}})^2}\sum_{i=1}^{\infty}y_{\alpha+1,i}^{2k-1}q^{i-4}\frac{E_5^{6i-2}}{E_1^{6i-4}}.
	\end{align*}
	Applying the $H$-operator to both sides of the above equation, we obtain
	\begin{align*}
		\sum_{n=0}^{\infty}RD_{(5^{2k-1},5^{2k-1})}^2\big(5^{2k+\alpha-1}(&5n+3)+\frac{9\cdot5^{2k+\alpha-1}-2\cdot5^{2k-1}+1}{12}\big)q^{5n} \nonumber\\=&\frac{1}{({E_{5^{2k-\alpha-1}}})^2}\sum_{i=1}^{\infty}y_{\alpha+1,i}^{2k-1}H\big(q^{i-4}\frac{E_5^{6i-2}}{E_1^{6i-4}}\big)\\
		=&\frac{1}{({E_{5^{2k-\alpha-1}}})^2}\sum_{i=1}^{\infty}y_{\alpha+1,i}^{2k-1}\sum_{j=1}^{5i-3} a_{i,j}\,
		q^{5j-5}\frac{E_{25}^{6j-2}}{E_{5}^{6j-4}}\\
		=&\frac{1}{({E_{5^{2k-\alpha-1}}})^2}\sum_{j=1}^{\infty}\big(\sum_{i=1}^{\infty}y_{\alpha+1,i}^{2k-1}a_{i,j}\big)q^{5j-5}\frac{E_{25}^{6j-2}}{E_{5}^{6j-4}}\\
		=&\frac{1}{({E_{5^{2k-\alpha-1}}})^2}\sum_{j=1}^{\infty}y_{\alpha+2,j}^{2k-1}q^{5j-5}\frac{E_{25}^{6j-2}}{E_{5}^{6j-4}}.
			\end{align*}
Replacing $q^5$ by $q$ implies that Theorem \ref{gf1} is true for $\alpha+1$, hence by induction, we complete the proof of Theorem \ref{gf1}.
\end{proof}

\begin{theorem}For each $n,\alpha\geq0$ and $k\geq1$,
	\begin{align}
		\sum_{n=0}^{\infty}RD_{(5^{2k-1},5^{2k-1})}^2\big(5^{4k+2\alpha-2}n+\frac{8\cdot5^{4k+2\alpha-2}+({5^{2k-1}-1})^2}{12}\big)q^n \nonumber\\
		=\sum_{i=1}^{\infty}z_{2\alpha+1,i}q^{i-1}\frac{E_5^{6i-2}}{E_1^{6i-2}},\label{3.2.1} \\
		\sum_{n=0}^{\infty}RD_{(5^{2k-1},5^{2k-1})}^2\big(5^{4k+2\alpha-1}n+\frac{4\cdot5^{4k+2\alpha-1}+({5^{2k-1}-1})^2}{12}\big)q^n \nonumber\\ =\sum_{i=1}^{\infty}z_{2\alpha+2,i}q^{i-1}\frac{E_5^{6i-4}}{E_1^{6i-4}}. \label{3.2.2}
		\end{align}
			where the coefficient entries are defined as follows:
		\begin{align*}
			Z_1=&Y_{2k-1}^{2k-1}\cdot A
		\end{align*}
		and for any $\alpha\geq1$, 
		\begin{align}\label{zdef}
				Z_{\alpha+1}=&\begin{cases}
				Z_\alpha\cdot A,& \text{if $\alpha$ is even},\\
				Z_\alpha\cdot B,& \text{if $\alpha$ is odd}
			\end{cases}
		\end{align}
	where,	\begin{align*}
		A:=&(a_{i,j})=m_{6i-4,i+j-1}\\
		B:=&(b_{i,j})=m_{6i-2,i+j-1}		
	\end{align*}
\end{theorem}
\begin{proof}
 For case, $\alpha=2k-2$ in Theorem \ref{gf1} gives
	\begin{align*}
		\sum_{n=0}^{\infty}RD_{(5^{2k-1},5^{2k-1})}^2\big(5^{4k-3}n+\frac{9\cdot5^{4k-3}-2\cdot5^{2k-1}+1}{12}\big)q^{n-3}\\ =\sum_{i=1}^{\infty}y_{2k-1,i}^{2k-1}q^{i-4}\frac{E_5^{6i-4}}{E_1^{6i-4}}.\\
	\end{align*}
	Applying $H$-operator in above identity, we have
	\begin{align*}
		\sum_{n=0}^{\infty}RD_{(5^{2k-1},5^{2k-1})}^2\big(5^{4k-3}(5&n+3)+\frac{9\cdot5^{4k-3}-2\cdot5^{2k-1}+1}{12}\big)q^{5n}\\ =&\sum_{i=1}^{\infty}y_{2k-1,i}^{2k-1}H\big(q^{i-4}\frac{E_5^{6i-4}}{E_1^{6i-4}}\big)\\
		=&\sum_{i=1}^{\infty}y_{2k-1,i}^{2k-1}\sum_{j=1}^{5i-3} a_{i,j}\, 
		q^{5j-5}{(\frac{E_{25}}{E_{5}})}^{6j-2}\\
		=&\sum_{j=1}^{\infty}\big(\sum_{i=1}^{\infty}y_{2k-1,i}^{2k-1}a_{i,j}\big)	q^{5j-5}{(\frac{E_{25}}{E_{5}})}^{6j-2}\\
		=&\sum_{j=1}^{\infty}z_{1,j}q^{5j-5}{(\frac{E_{25}}{E_{5}})}^{6j-2}.
		\end{align*}
	Finally, on	replacing $q^5$ by $q$ we have
		\begin{align*}
				\sum_{n=0}^{\infty}RD_{(5^{2k-1},5^{2k-1})}^2\big(5^{4k-2}n+\frac{45\cdot5^{4k-3}-2\cdot5^{2k-1}+1}{12}\big)q^{n}\\	=\sum_{j=1}^{\infty}z_{1,j}q^{j-1}{(\frac{E_{5}}{E_{1}})}^{6j-2}
		\end{align*}
		which is the $\alpha=0$ case of (\ref{3.2.1}).\\
		Now, suppose (\ref{3.2.1}) is true for some $\alpha\geq0$, we have
		\begin{align*}
			\sum_{n=0}^{\infty}RD_{(5^{2k-1},5^{2k-1})}^2\big(5^{4k+2\alpha-2}n+\frac{8\cdot5^{4k+2\alpha-2}+({5^{2k-1}-1})^2}{12}\big)q^{n-1} \\
			=\sum_{i=1}^{\infty}z_{2\alpha+1,i}q^{i-2}\frac{E_5^{6i-2}}{E_1^{6i-2}}.
			\end{align*}
				Applying $H$-operator to both sides, we have
				\begin{align*}
					\sum_{n=0}^{\infty}RD_{(5^{2k-1},5^{2k-1})}^2\big(5^{4k+2\alpha-2}(5n+1)+&\frac{8\cdot5^{4k+2\alpha-2}+({5^{2k-1}-1})^2}{12}\big)q^{5n} \\
					=&\sum_{i=1}^{\infty}z_{2\alpha+1,i}H\big(q^{i-2}\frac{E_5^{6i-2}}{E_1^{6i-2}}\big)\\
					=&\sum_{i=1}^{\infty}z_{2\alpha+1,i}\sum_{j=1}^{5i-1} b_{i,j}\,
					q^{5j-5}{(\frac{E_{25}}{E_{5}})}^{6j-4}\\
					=&\sum_{j=1}^{\infty}\big(\sum_{i=1}^{\infty}z_{2\alpha+1,i}b_{i,j}\big)q^{5j-5}{(\frac{E_{25}}{E_{5}})}^{6j-4}\\
					=&\sum_{j=1}^{\infty}z_{2\alpha+2,j}q^{5j-5}{(\frac{E_{25}}{E_{5}})}^{6j-4}\\
				\end{align*}
					Replacing $q^5$ by $q$ gives 
					\begin{align*}
							\sum_{n=0}^{\infty}RD_{(5^{2k-1},5^{2k-1})}^2\big(5^{4k+2\alpha-1}n+\frac{4\cdot5^{4k+2\alpha-1}+({5^{2k-1}-1})^2}{12}\big)q^n \\ =\sum_{i=1}^{\infty}z_{2\alpha+2,i}q^{i-1}\frac{E_5^{6i-4}}{E_1^{6i-4}}
					\end{align*}
	which is (\ref{3.2.2}).\\
	Suppose that (\ref{3.2.2}) is true for some $\alpha\geq0$
	\begin{align*}
		\sum_{n=0}^{\infty}RD_{(5^{2k-1},5^{2k-1})}^2\big(5^{4k+2\alpha-1}n+\frac{4\cdot5^{4k+2\alpha-1}+({5^{2k-1}-1})^2}{12}\big)q^{n-3} \\ =\sum_{i=1}^{\infty}z_{2\alpha+2,i}q^{i-4}\frac{E_5^{6i-4}}{E_1^{6i-4}}. 
	\end{align*}
	Applying $H$-operator gives
	\begin{align*}
		\sum_{n=0}^{\infty}RD_{(5^{2k-1},5^{2k-1})}^2\big(5^{4k+2\alpha-1}(5n+3)+&\frac{4\cdot5^{4k+2\alpha-1}+({5^{2k-1}-1})^2}{12}\big)q^{5n} \\ =&\sum_{i=1}^{\infty}z_{2\alpha+2,i}H\big(q^{i-4}\frac{E_5^{6i-4}}{E_1^{6i-4}}\big)\\
		=&\sum_{i=1}^{\infty}z_{2\alpha+2,i}\sum_{j=1}^{5i-3} a_{i,j}\,
		q^{5j-5}{(\frac{E_{25}}{E_{5}})}^{6j-2}\\
		=&\sum_{j=1}^{\infty}\big(\sum_{i=1}^{\infty}z_{2\alpha+2,i} a_{i,j}\big)q^{5j-5}{(\frac{E_{25}}{E_{5}})}^{6j-2}\\
		=&\sum_{j=1}^{\infty}z_{2\alpha+3,j}q^{5j-5}{(\frac{E_{25}}{E_{5}})}^{6j-2}\\
	\end{align*}
Replacing $q^5$ by $q$ 
\begin{align*}
	\sum_{n=0}^{\infty}RD_{(5^{2k-1},5^{2k-1})}^2\big(5^{4k+2\alpha}n+\frac{8\cdot5^{4k+2\alpha}+({5^{2k-1}-1})^2}{12}\big)q^n \\
	=\sum_{i=1}^{\infty}z_{2\alpha+3,i}q^{i-1}\frac{E_5^{6i-2}}{E_1^{6i-2}}
\end{align*}
which is the $\alpha+1$ case of (\ref{3.2.1}).\\This completes the proof of Theorem 3.2 by induction.
\end{proof}
\section{Proof of Theorems \ref{thm:main}-\ref{t2} }
We are now in a  position to give proof of Theorems 1.1 to 1.3. We begin by focusing on Theorem 1.1

\begin{proof}[\textbf{Proof of Theorem \ref{thm:main}}]
	Fix $k\ge1$ and let $\rho\in\{0,1,\dots,p-1\}$ be the residue when
	$\rho\equiv-k\pmod p$. Put $j=(k+\rho)/p\in\Z_{\ge1}$. From (\ref{Gf})
	\begin{align}
		\sum_{n=0}^{\infty}\RD{k}{n}q^n
		=&\left(\frac{E_lE_t}{E_{lt}}\right)^{k}\frac1{E_1^{k}}\nonumber\\
		=&\left(\frac{E_lE_t}{E_{lt}}\right)^{k}\frac{E_1^{\rho}}{E_1^{pj}}.\nonumber\\
			\equiv& G(q)\,E_1^{\rho}\pmod p,
		\qquad
		G(q):=\left(\frac{E_lE_t}{E_{lt}}\right)^{k}E_p^{-j}.\label{rho}
	\end{align}
	Since $E_l,E_t,E_{lt}$ and $E_p$ are series in $q^p$. Therefore, the terms of exponent congruent to r modulo p on the left side of (\ref{rho}) are congruent to $G(q)$ times the terms of $E_1^{\rho}$ of exponent congruent to $ r$ modulo $p$. Also, we note $G$ has only integeral coefficients only.\\
	\textbf{Case (1)} : For $\rho=0$ or $k\equiv0\pmod p$\\
	Setting $k=pm$, we have $j=m$ and $E_1^{\rho}=1$, whose only term has exponent $0$.\\
	Now if $1\le r\le p-1$, we have no term of exponent congruent to $r$ modulo $p$, so $\RD{k}{pn+r}\equiv0\pmod p$.\\
	Define
	$F:=\bigl(E_lE_t/(E_1E_{lt})\bigr)^{m}=\sum_{n\ge0}\RD{m}{n}q^n$.
	\begin{align*}
		\sum_{n\ge0}\RD{pm}{n}q^n=F^p\equiv F(q^p)=\sum_{n\ge0}\RD{m}{n}q^{pn}\pmod p .
	\end{align*}
	On comparing the coefficients of $q^{pn}$ we obtain 
	\begin{align*}
		\RD{pm}{pn}\equiv\RD{m}{n}\pmod p.
	\end{align*}
	\textbf{Case (2) }: For $\rho=1$ or $k\equiv-1\pmod p$\\
Employing Lemma \ref{pentagonal} in (\ref{rho}), we have
	\begin{align*}
		\sum_{n=0}^{\infty}\RD{k}{n}q^n
		\equiv&&G(q)\sum_{k \in \mathbb{Z}} q^{\frac{k(3k+1)}{2}} \pmod{p},
		\qquad
		G(q):=\left(\frac{E_lE_t}{E_{lt}}\right)^{k}E_p^{-j}.
	\end{align*}
	Now, we are interested in comparing the exponents of $q^{pn+r}$ and $q^{k(3k+1)/2}$. That
	is, we want to know if it is ever possible to have
	\begin{align*}
		pn+r = \frac{k(3k+1)}{2},
	\end{align*}
	for specific $n$ and $k$. Note that, if this is possible, then we must have
	\begin{align*}
		r \equiv \frac{k(3k+1)}{2} \pmod{p}
	\end{align*}
	Completing the square, this means we would have
	\begin{align*}
		24r + 1 \equiv (6k + 1)^2 \pmod{p}.
	\end{align*}
	However, this is never possible because we have assumed that 24r+1 is a quadratic
	nonresidue modulo $p$. This completes the proof of Case (2).\\
	\textbf{Case (3)} : For $\rho=3$ or $k\equiv-3\pmod p$\\
Employing Lemma \ref{jacobi} in (\ref{rho}), we have
	\begin{align*}
		\sum_{n=0}^{\infty}\RD{k}{n}q^n
		\equiv&G(q) \sum_{k=0}^{\infty} q^{\frac{k(k+1)}{2}} \pmod{p},
		\qquad
		G(q):=\left(\frac{E_lE_t}{E_{lt}}\right)^{k}E_p^{-j}.
	\end{align*}
	Now, we are interested in comparing the exponents of $q^{pn+r}$ and $q^{\frac{k(k+1)}{2}}$. That
	is, we want to know if it is ever possible to have
	\begin{align*}
		pn+r = \frac{k(k+1)}{2},
	\end{align*}
	for specific $n$ and $k$. Note that, if this is possible, then we must have
	\begin{align*}
		r \equiv \frac{k(k+1)}{2} \pmod{p}
	\end{align*}
	Completing the square, this means we would have
	\begin{align*}
		8r + 1 \equiv (2k + 1)^2 \pmod{p}.
	\end{align*}
	However, this is never possible because we have assumed that 8r+1 is a quadratic
	nonresidue modulo $p$. This completes the proof of Case (3).\\
	\textbf{Case (4)} : For $\rho=4$ or $p\equiv2\pmod{3}$ and $6r+1\equiv0\pmod p$\\
Employing Lemma \ref{pentagonal} and \ref{jacobi} in (\ref{rho}), we have
	\begin{align*}
		\sum_{n=0}^{\infty}\RD{k}{n}q^n
		\equiv& G(q)\,E_1^4\pmod p,
		\qquad
		G(q):=\left(\frac{E_lE_t}{E_{lt}}\right)^{k}E_p^{-j}.\\
		=&\sum_{a\in\Z}\ \sum_{b\ge0}(-1)^{a+b}(2b+1)\,
		q^{\,a(3a+1)/2+b(b+1)/2}.
	\end{align*}
	Now, we are interested in comparing the exponents of $q^{pn+r}$ and $q^{\,a(3a+1)/2+b(b+1)/2}$. That
	is, we want to know if it is ever possible to have
	\begin{align*}
		pn+r = \frac{a(3a+1)+b(b+1)}{2},
	\end{align*}
	for specific $n$ and $k$. Note that, if this is possible, then we must have
	\begin{align*}
		r \equiv \frac{a(3a+1)+b(b+1)}{2}\pmod{p}
	\end{align*}
	
Specifically, set $x=6a+1$ and $y=2b+1$, and also
	$24\cdot\frac{a(3a+1)}{2}=x^2-1$ and $24\cdot\frac{b(b+1)}{2}=3y^2-3$, we have
	\begin{align*}
		24r + 4 \equiv x^2+3y^2 \pmod{p}.
	\end{align*}
	Suppose $6r+1\equiv0\pmod p$, then
	$24r+4=4(6r+1)\equiv0\pmod p$, we have
	\[
	x^2+3y^2\equiv0\pmod p .
	\]
	Assume, $p\nmid y$, implies $y$ has an inverse $y'$ modulo $p$.
	Multiplying above identity by $y'^2$ we obtain
	\begin{align*}
		({xy'})^2\equiv{-3}\pmod p.
	\end{align*}
A contradiction since $p\equiv2\pmod3$, we have $\leg{-3}{p}=\leg{p}{3}=-1$.
	Hence $p\mid y=2b+1$, so the term $(-1)^{a+b}(2b+1)q^{a(3a+1)/2+b(b+1)/2}$ has coefficient
	divisible by $p$. Thus every term of $E_1^4$ with exponent congruent to $r$ modulo $p$ has
	coefficient divisible by $p$.\\
	This completes the proof of Theorem~\ref{thm:main}.
\end{proof}
\noindent Guaranteed by Theorem~\ref{thm:main} in the following table, we provide a complete classification for small primes and $r\geq1$, satisfying the congruence
$\RD{k}{pn+r}\equiv0\pmod p$. 
\begin{table}[H]
	\centering
	\renewcommand{\arraystretch}{1.3}
	\begin{tabular}{c c l l}
		\hline
		$p$ & $k \bmod p$ & Vanishing classes $r$ & Congruence\\
		\hline
		5  & $0$  & $1,2,3,4$            & $\mathrm{RD}_{l,t}^{5j}(5n+r)\equiv0\pmod 5$\\
		5  & $4$  & $3,4$                & $\mathrm{RD}_{l,t}^{5j+4}(5n+r)\equiv0\pmod 5$\\
		5  & $2$  & $2,3,4$              & $\mathrm{RD}_{l,t}^{5j+2}(5n+r)\equiv0\pmod 5$\\
		5  & $1$  & $4$                  & $\mathrm{RD}_{l,t}^{5j+1}(5n+4)\equiv0\pmod 5$\\
		\hline
		7  & $0$  & $1,\dots,6$          & $\mathrm{RD}_{l,t}^{7j}(7n+r)\equiv0\pmod 7$\\
		7  & $6$  & $3,4,6$              & $\mathrm{RD}_{l,t}^{7j+6}(7n+r)\equiv0\pmod 7$\\
		7  & $4$  & $2,4,5,6$            & $\mathrm{RD}_{l,t}^{7j+4}(7n+r)\equiv0\pmod 7$\\
		\hline
		11 & $0$  & $1,\dots,10$         & $\mathrm{RD}_{l,t}^{11j}(11n+r)\equiv0\pmod{11}$\\
		11 & $10$ & $3,6,8,9,10$         & $\mathrm{RD}_{l,t}^{11j+10}(11n+r)\equiv0\pmod{11}$\\
		11 & $8$  & $2,4,5,7,8,9$        & $\mathrm{RD}_{l,t}^{11j+8}(11n+r)\equiv0\pmod{11}$\\
		11 & $7$  & $9$                  & $\mathrm{RD}_{l,t}^{11j+7}(11n+9)\equiv0\pmod{11}$\\
		\hline
		13 & $0$  & $1,\dots,12$         & $\mathrm{RD}_{l,t}^{13j}(13n+r)\equiv0\pmod{13}$\\
		13 & $12$ & $3,4,6,8,10,11$      & $\mathrm{RD}_{l,t}^{13j+12}(13n+r)\equiv0\pmod{13}$\\
		13 & $10$ & $4,5,7,8,9,11,12$    & $\mathrm{RD}_{l,t}^{13j+10}(13n+r)\equiv0\pmod{13}$\\
		\hline
	\end{tabular}
\end{table}
We now require the complementary lemmas to comprimise the prove of Theorem \ref{T1} and \ref{t2}.
For any integer $n$, let $v(n)$ denote the highest power of $5$ that divides $n$ and $v(0)=\infty$, we have the following lemmas.

\begin{lemma}[{\cite[Eq.~(6.5.1)]{power}}]\label{4.1}
	For any $i\geq 1$ and $j\geq 1$,
	\begin{equation}\label{2.10}
		v(m_{i,j}) \geq \left[\frac{5j-i-1}{2}\right].
	\end{equation}
\end{lemma}
\begin{lemma}\label{lem4.2}
	For any $k\geq 1$, $0\leq\alpha\leq 2k-2$ and $j\geq 1$,
	\begin{equation}\label{4.3}
		v\big(y^{(2k-1)}_{\alpha+1,j}\big)\geq k+\left[\frac{\alpha}{2}\right]+\left[\frac{5j-5}{2}\right] ,
	\end{equation}
	where $[x]$ denotes the greatest integer less than or equal to $x$.
\end{lemma}

\begin{proof}
	It follows from \cite[Eq. (4.5)]{tang} that
	\[
	v\big(y^{(2k-1)}_{1,j}\big)=v(x_{2k-1,j})\geq k+\left[\frac{5j-5}{2}\right].
	\]
	This implies that \eqref{4.3} holds for $\alpha=0$.
	
	Now we assume that \eqref{4.3} is true for some $\alpha$ with $0\leq\alpha\leq 2k-3$.
	With the help of \eqref{Ydef}, we obtain that
	\begin{align*}
		v\big(y^{(2k-1)}_{\alpha+2,j}\big)
		&\geq \min_{i\geq 1}\Big\{v\big(y^{(2k-1)}_{\alpha+1,i}\big)+v(a_{i,j})\Big\}\\
		&\geq \min_{i\geq 1}\left\{k+\left[\frac{\alpha}{2}\right]+\left[\frac{5i-5}{2}\right]
		+\left[\frac{5j-i-2}{2}\right]\right\}\\
		&\geq k+\left[\frac{\alpha}{2}\right]+1+\left[\frac{5j-5}{2}\right]\\
		&\geq k+\left[\frac{\alpha+1}{2}\right]+\left[\frac{5j-5}{2}\right],
	\end{align*}
	This implies \eqref{4.3} holds for $\alpha+1$. We therefore complete the proof by induction.
\end{proof}

\begin{lemma}
		For any $k\geq 1$, $\alpha\geq 0$ and $j\geq 1$,
	\begin{align}
		v\big(z_{2\alpha+1,j}\big) &\geq 2k+\alpha+\left[\frac{5j-5}{2}\right] \label{2.11}\\
	v\big(z_{2\alpha+2,j}\big) &\geq 2k+\alpha+\left[\frac{5j-5}{2}\right] \label{2.12}
	\end{align}
	
\end{lemma}
\begin{proof}
From (\ref{zdef}), we have 
\begin{align*}
	v(z_{1,j})=&v(\sum_{i=1}^{\infty}y_{2k-1,i}^{2k-1}\cdot a_{i,j})\\
	\leq&2k+\left[\frac{5j-5}{2}\right]
\end{align*}
This implies that \eqref{2.11} holds for $\alpha=0$.

Now we assume that \eqref{2.11} is true for some $\alpha\geq0$. With the help of (\ref{zdef}), we obtain that
\begin{align*}
	v\big(z_{2\alpha+2,j}\big)
	&\geq \min_{i\geq 1}\left\{v\big(z_{2\alpha+1,i}\big)+v(b_{i,j})\right\}\\
	&\geq \min_{i\geq 1}\left\{2k+\alpha+\left[\frac{5i-5}{2}\right]+\left[\frac{5j-i-4}{2}\right]\right\}\\
	&\geq 2k+\alpha+\left[\frac{5j-5}{2}\right],
\end{align*}
which is \eqref{2.12}.

Now suppose \eqref{2.12} holds for some $\alpha\geq 0$. Then,
\begin{align*}
	v\big(z_{2\alpha+3,j}\big)
	&\geq \min_{i\geq 1}\left\{v\big(z_{2\alpha+2,i}\big)+v(a_{i,j})\right\}\\
	&\geq \min_{i\geq 1}\left\{2k+\alpha+\left[\frac{5i-5}{2}\right]+\left[\frac{5j-i-2}{2}\right]\right\}\\
	&\geq 2k+\alpha+1+\left[\frac{5j-5}{2}\right].
\end{align*}
This implies \eqref{2.11} holds for $\alpha+1$. We therefore complete the proof by induction.
\end{proof}
\begin{proof}[\textbf{Proof of Theorem \ref{T1}}]
	It follows from emma \ref{lem4.2} that
	\begin{align}
	v\big(y^{(2k-1)}_{\alpha+1,j}\big)\geq
	\begin{cases}
		k+\left[\frac{\alpha}{2}\right] & \text{if } j=1,\\[2pt]
		k+\left[\frac{\alpha}{2}\right]+2 & \text{if } j\geq 2.
	\end{cases}
	\end{align}
	The congruence (\ref{t11}) follows immediately from above inequality and Theorem \ref{gf1}.
	
In particular taking, modulo $5^{	k+\left[\frac{\alpha}{2}\right]+1}$,
	\begin{align}
\sum_{n=0}^{\infty}RD_{(5^{2k-1},5^{2k-1})}^2\big(5^{2k+\alpha-1}n+\frac{9\cdot5^{2k+\alpha-1}-2\cdot5^{2k-1}+1}{12}\big)q^n \nonumber\\ \equiv\frac{1}{({E_{5^{2k-\alpha-1}}})^2}y_{\alpha+1,1}^{2k-1}{E_5^{3}}{E_1^{3}}.\label{extra11}
	\end{align}
	With the help of \eqref{r(q)}, we obtain that
	\begin{equation}\label{f1^31}
		E_1^{3}=E_{25}^{3}\big(R(q^{5})^{-3}-3qR(q^{5})^{-2}+5q^{3}-3q^{5}R(q^{5})^{2}-q^{6}R(q^{5})^{3}\big).
	\end{equation}
	The congruences \eqref{t12} and \eqref{t14} follow by utilizing (\ref{extra11}) and \eqref{f1^31}.
	This completes the proof of Theorem \ref{T1}.
\end{proof}
\begin{proof}[\textbf{Proof of Theorem \ref{t2}}]
It follows from (\ref{2.11}) and (\ref{2.12}) that
	\begin{align}
		\nu\big(z_{2\alpha+1,j}\big) &\geq \label{vrqz1}
		\begin{cases}
			2k+\alpha & \text{if } j=1,\\
			2k+\alpha+2 & \text{if } j\geq 2,
		\end{cases} \\[4pt]
		\nu\big(z_{2\alpha+2,j}\big) &\geq \label{vrqz2}
		\begin{cases}
			2k+\alpha & \text{if } j=1,\\
			2k+\alpha+2 & \text{if } j\geq 2.
		\end{cases}
	\end{align}
	The congruence (\ref{t21}) follows immediately from (\ref{3.2.1}) and (\ref{vrqz1}). And congruence (\ref{t22}) follows from (\ref{3.2.2}) and (\ref{vrqz2}).\\
	Similarly taking modulo $5^{2k+\alpha+1}$,
	\begin{align}
		\sum_{n=0}^{\infty}RD_{(5^{2k-1},5^{2k-1})}^2\big(5^{4k+2\alpha-2}n+\frac{8\cdot5^{4k+2\alpha-2}+({5^{2k-1}-1})^2}{12}\big)q^n \nonumber\\
		\equiv z_{2\alpha+1,1}{E_5^3}{E_1} \label{extra1}
	\end{align}
	Thanks to (\ref{r(q)}), we have
	\begin{equation}\label{2.24}
		E_1=E_{25}\big(R(q^{5})^{-1}-q-q^{2}R(q^{5})\big),
	\end{equation}
	where
	\[
	R(q)=\frac{(q;q^{5})_{\infty}(q^{4};q^{5})_{\infty}}{(q^{2};q^{5})_{\infty}(q^{3};q^{5})_{\infty}}.
	\]
	The congruences (\ref{t23}) and (\ref{t24}) follow by combining (\ref{extra1}) and \eqref{2.24}.
	
	From (\ref{vrqz2}) and \eqref{3.2.2}, and taking modulo $5^{2k+\alpha+1}$, we obtain
	\begin{align}
		\sum_{n=0}^{\infty}RD_{(5^{2k-1},5^{2k-1})}^2\big(5^{4k+2\alpha-1}n+\frac{4\cdot5^{4k+2\alpha-1}+({5^{2k-1}-1})^2}{12}\big)q^n\nonumber \\
		\equiv z_{2\alpha+2,1}{E_5}{E_1^3}.\label{extra2}
	\end{align}
	On employing \eqref{2.24}, we have
	\begin{equation}\label{f1^3}
		E_1^{3}=E_{25}^{3}\big(R(q^{5})^{-3}-3qR(q^{5})^{-2}+5q^{3}-3q^{5}R(q^{5})^{2}-q^{6}R(q^{5})^{3}\big).
	\end{equation}
	The congruences \eqref{t25} and \eqref{t26} follow by utilizing \eqref{extra2} and \eqref{f1^3}.
	This completes the proof of Theorem \ref{t2}.
\end{proof}

\noindent\textsuperscript{1,2}Ramanujan School of Mathematical Sciences,\\Department of Mathematics,\\ Pondicherry University,\\ Puducherry - 605014, India. \\

\noindent Email: \texttt{anjelinvallialil@pondiuni.ac.in} \\
Email: \texttt{dr.fathima.sn@pondiuni.ac.in}  (\Letter)

\begin{thebibliography}{99}
	\bibitem{power}
	Hirschhorn, M. D, The power of q: a personal journey, \textit{Springer} (2017).
\bibitem{Hunt}
Hirschhorn, M. D, Hunt, D. C, A simple proof of the Ramanujan conjecture for powers of 5 (1981).
\bibitem{Anjelin}
Johnson, A. M.,  Fathima, S. N., On Congruences modulo powers of 3 for Simultaneously $l$-Regular and $t$-Distinct Partitions in $k$-Colors, arXiv preprint arXiv:2607.11499 (2026).
\bibitem{keith}
Keith, William J., Partitions into Parts Simultaneously Regular, Distinct, And/or Flat, {\it Springer Proceedings in Mathematics  Statistics, Springer International Publishing} (2017).
\bibitem{keith2}
Keith, W. J., A bijection for partitions simultaneously $ s $-regular and $ t $-distinct, arXiv preprint arXiv:2207.13840 (2022).
\bibitem{Ahmia 2}
Nadji, M. L., Ahmia, M., Partitions into parts simultaneously regular and distinct, {\it Journal of the Ramanujan Mathematical Society}, {40}(2), (2025).
\bibitem{tang}
Tang, D, Congruences modulo powers of 5 for k-colored partitions, \textit{ Journal of Number Theory}, \textbf{187}, 198-214 (2018).

\end{thebibliography}
\end{document}